\documentclass[a4paper, 10pt]{article}
\usepackage{anysize}
\marginsize{3,5cm}{2,5cm}{2,5cm}{2,5cm}
\usepackage{amssymb}
\usepackage{amsmath,amsbsy,amssymb,amscd}
\usepackage{t1enc}
\usepackage[cp1250]{inputenc}
\usepackage[british]{babel}
\usepackage[all]{xy}
\usepackage{color}
\usepackage{amsfonts}
\usepackage{latexsym}
\usepackage{amsthm}
\usepackage{mathrsfs}

\usepackage{mathrsfs} 

\DeclareMathAlphabet{\mathpzc}{OT1}{pzc}{L}{it} 

\newtheorem{definition}{Definition}[section]

\newtheorem{proposition}[definition]{Proposition}
\newtheorem{theorem}[definition]{Theorem}
\newtheorem{corollary}[definition]{Corollary}
\newtheorem{remark}[definition]{Remark}
\newtheorem{lemma}[definition]{Lemma}

\def\geq{\geqslant}
\def\leq{\leqslant}
\def\R{\mathbb{R}}
\def\T{\mathbb{T}}

\def\Z{\mathbb{Z}}
\def\N{\mathbb{N}}

\def\Q{\mathbb{Q}}
\newcommand{\bea}{\begin{eqnarray}}
  \newcommand{\eea}{\end{eqnarray}}
  \newcommand{\beab}{\begin{eqnarray*}}
  \newcommand{\eeab}{\end{eqnarray*}}

  \newcommand{\be}{\begin{equation}}
  \newcommand{\ee}{\end{equation}}
\title{Ratner's property for special flows over irrational rotations under functions of bounded variation. II}
\author{Adam Kanigowski}
\begin{document}
\baselineskip=14pt \maketitle
\begin{abstract} We consider special flows over the rotation on the circle by an irrational $\alpha$ under roof functions of bounded variation. The roof functions, in the Lebesgue decomposition, are assumed to have a continuous singular part coming from a quasi-similar Cantor set (including the Devil's staircase case). Moreover, a finite number of discontinuities is allowed. Assuming that $\alpha$ has bounded partial quotients, we prove that all such flows are 
weakly mixing and enjoy weak Ratner's property. Moreover, we provide a sufficient condition for the roof function to obtain a stability of the cocycle Ratner's property for the resulting special flow. 

\end{abstract}

\section{Introduction}  \indent The paper is a continuation of \cite{Kan}, that is, we will continue the study of Ratner's property (originally, H-property in \cite{Rat}) in the class of special flows $\mathcal{T}^f:=(T_t^f)_{t\in\R}$ determined by the rotation $Tx=x+\alpha$, on the additive circle $\T$, by an irrational number $\alpha$ with bounded partial quotients and $f$ a function of bounded variation. In order to study Ratner's property in \cite{Kan}, we used the Lebesgue decomposition 
\begin{equation}\label{2}f=f_j+f_a+f_s +\tilde{S}\{\cdot\},\end{equation}
where $f_j$ is the jump function (with countably many jumps $d_i$), $f_a$ is absolutely continuous on $\T$, $f_s$ is singular and continuous on $[0,1]$, $\tilde{S}:=\int_\T f'd\lambda$ ($\lambda$ denotes Lebesgue measure on $\T$) and $\{t\}$ stands for the fractional part of $t\in \R$. In \cite{Kan}, we assumed that 
\begin{equation}\label{ck} f_s=0\;\; \text{     and     }\;\; \tilde{S}\neq 0
\end{equation}
 (in that case $\tilde{S}=\sum_{i=1}^{+\infty}d_i$). We recall that some particular cases of (\ref{ck}) have already been studied in \cite{Fr-Lem,Fr-Lem-Les}. In \cite{Kan}, under the assumptions (\ref{ck}),  we proved the weak mixing property of $\mathcal{T}^f$ and, under an additional condition concerning the rate of convergence of the series of jumps, so called weak Ratner's property (WR-property from now on)\footnote{WR-property has the same dynamical consequences as the original H-property of Ratner, see Section~3.}. We also proved that (\ref{ck}) together with a condition on the set of jumps yield a stability result for WR-property in the class of roof functions studied in \cite{Kan} (for sufficiently small bounded variation perturbations satisfying (\ref{ck})).\\
\indent The situation becomes more complicated when $f_s\neq 0$. We have already noticed in \cite{Kan} that, due to \cite{Iw-Lem-Ma, Vol}, it may happen that for an irrational rotation there exists $f=f_s$, $f_s(1)-f_s(0)\neq 0$ such that the corresponding special flow $\mathcal{T}^f$ is isomorphic to the suspension flow (in particular such a special flow is not even weakly mixing).\\
\indent In the present paper, we constantly assume $f_s\neq 0$ and focus on the following two problems:
\begin{enumerate}
    \item WR-property and weak mixing for functions of bounded variation with the singular part being of the Devil's staircase type.
	 \item The stability of WR-property for arbitrary (sufficiently small) bounded variation perturbations, answering a question of  M.\ Lemańczyk \cite{Lem}.
\end{enumerate}
We now pass to a description of the results.\\
\indent We deal with the case when $f_s$ is naturally defined by a Cantor set. More precisly, we consider quasi-similar  Cantor sets. Each such set is obtained in the following way. We are given two bounded sequences of natural numbers $(m_i)_{i\geq 1}$, $(k_i)_{i\geq 1}$ satisfying
$$ 2\leq m_i\leq \left[\frac{k_i}{2}\right]~\footnote{By $[x]$ we denote the integer part of $x\in\R$.}, \; i\geq 1.$$
First, divide $[0,1)$ into $k_1$ clopen subintervals of equal length and choose $m_1$ of them (including the first one and the last one) so that the closures (in $[0,1]$) of any two of them are disjoint. Then set $A_1\subset [0,1]$ to be the union of the closures of the selected $m_1$ intervals. In the next step, we divide each of the previously selected intervals into $k_2$ clopen subintervals of equal length and we select $m_2$ subintervals (including the first one and the last one in each previously chosen subinterval), with the same configuration choice in each selected subinterval, so that the closures of any two of them are disjoint. We then set $A_2$ to be the union of the closures of all selected subintervals of step two. Clearly $A_2\subset A_1$. Proceeding in the same manner, we obtain a sequence $A_1\supset A_2\supset...\supset A_n\supset...$ of closed subsets of $[0,1]$ and we define the corresponding {\em quasi-similar} Cantor set $\mathcal{C}$ as $\mathcal{C}:=\bigcap_{i=1}^{+\infty}A_i$. Of course, it has Lebesgue measure $0$. As in the classical case, there is a canonicaly associated to $\mathcal{C}$ continuous on [0,1] non-decreasing function $f=f(\mathcal{C}):[0,1]\to \R$, $f(0)=0$, $f(1)=1$ whose derivative $f'=0$ on $[0,1]\setminus \mathcal{C}$ (an alternative possibility is the non-increasing function $f$ with the above properties but $f(0)=1, f(1)=0$).\\
\indent In the paper, we consider roof functions $f$ of bounded variation and Lebesgue decomposition (\ref{2})
 $$f=f_j+f_a+f_s+\tilde{S}\{\cdot\},$$
where
\begin{enumerate}
    \item $f_j$ has finitely  many discontinuity points $\{\beta_1,...,\beta_k\}$ (with the corresponding set of jumps $\{d_1,...,d_k\}$),
	 \item $f_s=f(\mathcal{C})$ (we consider both possibilities $f_s(0)=0$ or $f_s(0)=1$)\\ 
 and either  
 	 \item  $\tilde{S}\neq 0$\\
\noindent{or}
	 \item $\tilde{S}=0$, $f_j=0$, that is, $f=f_s+f_a$ (notice that this includes the classical Devil's staircase). 
\end{enumerate}
Note that $\tilde{S}=S-d_0$, where $S$ is the sum of jumps of $f$ ($S:=\sum_{i=0}^kd_i$, where $d_0=f_s(1)-f_s(0)$). The main results of the paper are summarized in the following.

\begin{theorem}\label{wmrp} Assume that $Tx=x=\alpha$ is the rotation by an irrational $\alpha\in [0,1)$ having bounded partial quotients. Let $f$ be a bounded variation roof function, $f=f_a+f_j+f_s+\tilde{S}\{\cdot\}$, where $f_s=f(\mathcal{C})$ and $\mathcal{C}$ is a quasi-similar Cantor set. If either
\begin{enumerate}
\item $f_j$ has finitely many discontinuities, $\tilde{S}>0$ and $f_s$ is non-decreasing~\footnote{The same result is obtained for $\tilde{S}<0$ and $f_s$ non-increasing.}\\
or
\item $f=f_s+f_a$,
\end{enumerate}then the special flow $(T_t^f)_{t\in\R}$ is weakly mixing and has WR-property.
\end{theorem}
The main tool used to prove WR-property for special flows over rotations, so called {\em cocycle Ratner's property}, has been introduced in \cite{Fra-Lem} (we also used it in \cite{Kan}). We recall its definition now.
 Let $X$ be a compact metric (the metric is denoted by $d$) Abelian group with $\mathscr{B}$ the $\sigma$-algebra of Borel subsets of $X$ and  $\mu$ Haar measure on $X$. Let $Tx=x+x_0$ be an ergodic rotation on $(X,\mathscr{B},\mu)$. Assume that $f\in L^{\infty}(X,\mathscr{B},\mu)$, $f\geq c$ for some constant $c>0$. Assume $P\subset\R\setminus\{0\}$ is a compact set.
 \begin{definition}{\em \cite{Fra-Lem}\label{krp} One says that the special flow $(T_t^f)$ given by $T$ and $f$ has the {\em cocycle Ratner's property} (with the set $P$) if  for every $\epsilon>0$ and $N\in\N$ there exist $\kappa=\kappa(\epsilon)>0$, $\delta=\delta(\epsilon,N)>0$\footnote{\label{zale} If neccesary, we will write $\kappa_f(\epsilon), \delta_f(\epsilon,N)$, etc. to emphasize the dependence of these parameters on $f$.} and a subset $Z=Z(\epsilon,N)\in\mathscr{B}$ with $\mu(Z)>1-\epsilon$ such that if $x,x'\in Z$, $0<d(x,x')<\delta$, then there are $M=M(x,x')=M(\epsilon,N,x,x')\geq N$, $L=L(x,x')=L(\epsilon,N,x,x')\geq N$ such that $\frac{L}{M}\geq \kappa$ and there exists $p=p(\epsilon,N,x,x')\in P$ such that 
$$\frac{1}{L}\left|\{n\in\Z\cap[M,M+L] :\; |f^{(n)}(x)-f^{(n)}(x')-p|<\epsilon\}\right|>1-\epsilon.\footnote{ We denote $f^{(n)}(x)=f(x)+f(Tx)+...+f(T^{n-1}x)$ for $n\geq 1$.}.$$} 
\end{definition}
As proved in \cite{Fr-Le}, if $(T_t^f)_{t\in\R}$ satisfies the cocycle Ratner's property and is weakly mixing, then it enjoys WR-property. All special flows considered in \cite{Fr-Lem,Fra-Lem,Fr-Lem-Les,Kan} are weakly mixing and satisfy the cocycle Ratner's property. We note in passing that, assuming weak mixing one can ask whether WR-property is in fact equivalent to the cocycle Ratner's property, but we will not consider this problem here, postponing it to the forthcoming paper \cite{Kan2}.\\
\indent In the notes \cite{Lem}, M.\ Lemańczyk formulated a natural question whether the cocycle Ratner's property is stable under sufficiently small bounded variation perturbations. We will provide the positive answer to this question (see Section 5) assuming that:
\begin{itemize}
    \item The sets $Z=Z(\epsilon,N)$ in Definition \ref{krp} are equal to $X$ for each $\epsilon>0$, $N\in \N$.
	 \item For every $\epsilon>0, N\in\N$ and $x,x'\in X$ the (partial) functions $(x,x')\to M(\epsilon,N,x,x')d(x,x')$ defined for $x\neq x'$, $d(x,x')<\delta$, $\epsilon>0$, $N\in \N$ are bounded away from $0$ and $\infty$ (uniformly in $\epsilon, N, x,x'$).
	 \item For every $\epsilon>0$, $N\in \N$ and $x,x'$ as above, we have $|f^{(n)}(x)-f^{(n)}(x')-p|<\epsilon$ for each $n\in \Z\cap [M,M+L]$.\footnote{All examples from \cite{Fr-Lem,Fr-Lem-Les,Kan} enjoy both above properties.} 
 \end{itemize}
In Theorem \ref{sta} below, we prove that under the above three assumptions, whenever $(T_t^f)_{t\in\R}$ satisfies the cocycle Ratner's property, so does  $(T_t^{f+g})_{t\in\R}$, where $g$ has sufficiently small variation.\\
\indent Finally, we notice that using the methods from \cite{Fr-Lem}, the assertion of Theorem \ref{wmrp} in case 1. can be strenghtend to the mild mixing property of the corresponding special flows.\\

\indent The author would like to thank Professor Mariusz Lemańczyk for his patience and guidence.

\section{Basic notions}
We will use notation from \cite{Kan}. We denote by $\T$ the circle group $\R/\Z$ which will be identified with the interval $[0,1)$ with addition ${\rm mod}\, 1$. For a real number $t$ denote by $\|t\|$ its distance to the nearest integer number (note that $\|t\|=\|-t\|$ and $\|qt\|\leq q\|t\|$ for $q\in\N$). For an irrational $\alpha\in\T$ denote by $(q_n)_{n=0}^{+\infty}$ its sequence of denominators, that is, we have 
$$\frac{1}{2q_nq_{n+1}}<\left|\alpha-\frac{p_n}{q_n}\right|<\frac{1}{q_nq_{n+1}},$$
where $q_0=1,\; q_1=a_1,\;\; q_{n+1}=a_{n+1}q_n+q_{n-1}$, $p_0=0,\; p_1=1,\;\; p_{n+1}=a_{n+1}p_n+p_{n-1}$ and $[0;a_1,a_2,...]$ stands for the continued fraction expansion of $\alpha$. One says that $\alpha$ has {\em bounded partial quotients} if the sequence $(a_n)_{n=1}^{+\infty}$ is bounded. In this case, if we set $C:=\sup\{a_n:\; n\in\N\}+1$ then $q_{n+1}\leq Cq_n$ and
$$\frac{1}{2Cq_n}\leq\frac{1}{2q_{n+1}}<\|q_n\alpha\|<\frac{1}{q_{n+1}}<\frac{1}{q_n}$$
for each $n\in\N$. The following lemma is well-known.
\begin{lemma}\label{dist} Let $\alpha\in\T$ be irrational with bounded partial quotients. Then there exist positive constants $C_1,C_2$ such that for every $k\in\N$ the lengths of intervals $J_1,...,J_k$ arisen from the partition of $\T$ by $0,-\alpha,...,-(k-1)\alpha$ satisfy $\frac{C_2}{k}\leq|J_j|<\frac{C_1}{k}$ for each $j=1,...,k$. 
\end{lemma}
\begin{lemma}\label{q2} Let $\alpha$ be irrational with bounded partial quotients. Consider points $x_1,...,x_l\in \T$ such that there exists $q\in \Q$ so that for every $i,j\in 1,...l$, $x_i-x_j=\frac{r_{ij}}{q}$, with $|r_{ij}|\in \N\cap[0,q]$. If $C_2$ is a constant given by Lemma \ref{dist}, then for every $m\in \N$ the length of each interval arisen from the partition of $\T$ by the points $x_s+j\alpha$,  $s=1,...,l$ and $j=0,...,m-1$ is at least $\frac{C_2}{q^2m}$.
\end{lemma}
Proof: Consider any two points $x_i+t\alpha, x_j+s\alpha$ as above. Then
$$q\| (x_i+t\alpha)-(x_j+s\alpha)\|=q\|\frac{r_{ij}}{q}+ (t-s)\alpha\|\geq\|q(t-s)\alpha\|\geq \frac{C_2}{qm},$$ 
by Lemma \ref{dist} applied to $k=qm$. \hfill $\square$\\
\begin{remark}{ \em \label{dst} Take any $z_1,z_2\in \T$, $z_1-z_2=\frac{r}{q}$ for some $|r|<q$ ($r\in \Z$). Fix an interval $I\subset\T$ and assume that for some integers $t_0,t_1$, $t_0<t_1$, we have $z_1+t_0\alpha, z_2+t_1\alpha\in I$. Then 
\begin{equation}\label{100} t_1-t_0>\frac{C_2}{q^2|I|}.\end{equation}
Indeed, consider Lemma \ref{q2} in which $l=2$, $x_1=z_1+t_0\alpha$, $x_2=z_2+t_0\alpha$ and $m=t_1-t_0$. It follows that 
$$|I|\geq \|\frac{r}{q}+ (t_1-t_0)\alpha\|\geq \frac{C_2}{q^2(t_1-t_0)},$$ 
whence (\ref{100}) holds.} 
\end{remark}
\indent Assume that $T$ is an ergodic automorphism on $(X,\mathscr{B},\mu)$. We will always assume that $T$ is also aperiodic. A measurable function $f:X\to\R$ determines a cocycle $f^{(\cdot)}(\cdot):\Z\times X\to \R$ given by $f^{(m)}(x)=f(x)+f(Tx)+...+f(T^{m-1}x), m\geq 1,$ and
$$f^{(m+n)}(x)=f^{(m)}(x)+f^{(n)}(T^mx)$$
for each $m,n\in \Z$.
If $f:X\to\R$ is a strictly positive $L^1$ function, then by $\mathcal{T}^f=(T_t^f)_{t\in\R}$ we will mean the corresponding special flow under $f$ acting on $(X^f,\mathscr{B}^f,\mu^f)$, where $X^f:=\{(x,s)\in X\times\R:\; x\in X, 0\leq s<f(x)\}$, and $\mathscr{B}^f$ ($\mu^f$) is the restriction of $\mathscr{B}\otimes \mathscr{B}(\R)$ ($\mu\times \lambda$) to $X^f$. Under the action of the flow $\mathcal{T}^f$, each point in $X^f$ moves vertically with unit speed and we identify the point $(x,f(x))$ with $(Tx,0)$. More precisely, if $(x,s)\in X^f$ then
$$T_t^f(x,s)=(T^nx,s+t-f^{(n)}(x)),$$
where $n\in\Z$ is unique such that
$f^{(n)}(x)\leq s+t<f^{(n+1)}(x).$
 Assume additionally that $X$ is a metric space with a metric $d$. Then $X^f$ is a metrizable space with the metric
$$d_1((x,t),(y,s))=d(x,y)+|t-s|.$$
It is not hard to see that $\mathcal{T}^f$ satisfies the following ``almost continuity'' condition (\cite{Fra-Lem}): for every $\epsilon>0$ there exists $X(\epsilon)\in \mathscr{B}^f$, $\mu^f(X(\epsilon))>1-\epsilon$ such that for every $\epsilon'>0$ there exists $\epsilon_1>0$ such that 
$$d_1(T^f_u(x,t), T^f_{u'}(x,t))<\epsilon'$$
for each $(x,t)\in X(\epsilon)$ and $u,u'\in [-\epsilon_1,\epsilon_1]$.

\begin{proposition}\label{koks}{\em(Denjoy-Koksma inequality; see e.g. \cite{Kui-Nid})} If $f:\T\to\R$ is a function of bounded variation then
$$\left|\sum_{k=0}^{q_n-1}f(x+k\alpha)-q_n\int_\T f\,d\lambda\right|\leq {\rm Var}f,$$
for every $x\in \T$ and $n\in \N$.
\end{proposition}

\label{weak}\section{Weak Ratner's property}
In this section we recall  the notion of weak Ratner's property (WR-property) introduced in \cite{Fra-Lem} and we list results from \cite{Fra-Lem} needed in what follows.
Let $X$ be a $\sigma$-compact metric space with a metric $d$. Let $\mathscr{B}$ denote the  $\sigma$-algebra of Borel sets and $\mu$ a probability measure on $(X,\mathscr{B})$. Assume that $\mathcal{S}=(S_t)_{t\in\R}$ is a flow on $(X,\mathscr{B},\mu)$. Let $P\subset \R\setminus\{0\}$ be compact and $t_0\in \R\setminus\{0\}$. 

\begin{definition}{\em \cite{Fra-Lem}  The flow $(S_t)_{t\in\R}$ is said to have the property $R(t_0,P)$ if for every $\epsilon>0$ and $N\in\N$ there exist $\kappa=\kappa(\epsilon)$, $\delta=\delta(\epsilon,N)>0$ and a subset $Z=Z(\epsilon,N)\in\mathscr{B}$ with $\mu(Z)>1-\epsilon$ such that if $x,x'\in Z$, $x'$ is not in the orbit of $x$, and $d(x,x')<\delta$, then there are $M=M(x,x')\geq N$, $L=L(x,x')\geq N$ such that $\frac{L}{M}\geq \kappa$ and there exists $p=p(x,x')\in P$ such that 
$$\frac{1}{L}\big|\{n\in\Z\cap[M,M+L] :\; d(S_{nt_0}(x),S_{nt_0+p}(x'))<\epsilon\}\big|>1-\epsilon.$$
Moreover, we say that $(S_t)_{t\in\R}$ has  {\em weak Ratner's property}  or {\em WR-property} if the set of $s\in\R$ such that the flow $(S_t)_{t\in\R}$ has the $R(s,P)$-property is uncountable.}
\end{definition}

\begin{theorem}{\em \cite{Fra-Lem}} \label{join}
Let $(X,d)$ be a $\sigma$-compact metric space, $\mathscr{B}$ the $\sigma$-algebra of Borel subsets of $X$ and $\mu$ a probability Borel measure on $(X,\mathscr{B})$. Let $(S_t)_{t\in\R}$ be a weakly mixing flow on $(X,\mathscr{B},\mu)$ that satisfies the WR-property with a compact subset $P\subset\R\setminus\{0\}$. Assume that $(S_t)_{t\in\R}$ satisfies the ``almost continuity'' condition\footnote{ That is, $(S_t)_{t\in\R}$ satisfies ``almost continuity'' condition as the one for $\mathcal{T}^f$. The ``almost continuity'' condition which we gave for $\mathcal{T}^f$ can be naturally carried over to a more general situation.}. Let $(T_t)_{t\in\R}$ be an ergodic flow on $(Y,\mathscr{C},\nu)$  and let $\rho$ be an ergodic joining of $(S_t)_{t\in\R}$ and $(T_t)_{t\in\R}$. Then either $\rho=\mu\times \nu$ or $\rho$ is a finite extension of $\nu$.
\end{theorem}
We recall that each special flow $\mathcal{T}^f=(T_t^f)_{t\in\R}$ over an ergodic base $T$ is ergodic. In particular, the set $\{t\in\R; T_t^f\;\; \text{is not ergodic}\}$ is countable.
\begin{proposition}\label{condi}{\em \cite{Fra-Lem}} Let $X$ be a compact, metric (with a metric $d$) Abelian group with Haar measure $\mu$. Assume that $Tx=x+x_0$ is an ergodic rotation of $(X,\mathscr{B},\mu)$. Let $f\in L^{\infty}(X,\mathscr{B},\mu)$ be positive bounded away from zero. Assume that $\mathcal{T}^f=(T_t^f)_{t\in\R}$ satisfies the cocycle Ratner's property with a compact subset $P\subset \R\setminus\{0\}$. Suppose that $\gamma\in \R$ is a positive number such that the $\gamma$-time automorphism $T^f_\gamma:X^f\to X^f$ is ergodic. Then the special flow $\mathcal{T}^f$ has the $R(\gamma,P)$-property. In particular, $\mathcal{T}^f$ enjoys WR-property.
\end{proposition}
\subsection{Properties of bounded variation functions}\label{sec}
Let $T:\T\to\T$ be an irrational rotation by $\alpha$ with the sequence of denominators $(q_n)_{n=1}^{+\infty}$ and $f:\mathbb{T} \to \R$ a function of bounded variation. It follows by the Lebesgue decomposition, see e.g. \cite{Fre} and \cite{Kan}, that $f$ can be writen as
\begin{equation}\label{eq} f=f_a+f_j+f_s+\tilde{S}\{\cdot\},\end{equation}
where $f_j$ is the jump function on $\T$ with the set of jumps $\{d_i\}_{i=1}^{+\infty}$ at $\{\beta_i\}_{i=1}^{+\infty}$ respectively ($\beta_1=0$ and $d_1$ can be equal to zero), $f_s$ is singular, continuous on $[0,1]$ ($d_0:=f_s(1)-f_s(0)$) and $f_a$ is absolutely continuous on $\T$. Moreover, $\tilde{S}=\int_\T f'd\lambda$. Set $S:=\sum_{i=0}^{+\infty}d_i$ (then $\tilde{S}+d_0=S$) to be the sum of (all) jumps of $f$. Recall that the space ${\rm BV}(\T)$ of functions with bounded variation is a Banach space with the norm $\|f\|:={\rm Var}f+ \|f\|_{L^1}$.\\
\indent We will now focus on singular continuous functions which come from quasi-similar Cantor sets (see Introduction). By contruction, each such set is of the form $\mathcal{C}=\cap_{i=1}^{+\infty}A_i$ and for each $i\in \N$ the endpoints of the  chosen intervals in $A_i$ are always multiples of $\frac{1}{k_1...k_i}$, and there are $m_1...m_i$ of them. For each natural $n\in \N$ we define an absolutely continuous function $f_n$ such thatL: $f_n'=0$ on $[0,1]\setminus (\cup_{i=1}^n A_i)$ and $f_n$ is linear on each of the intervals of $A_n$, moreover $f_n(0)=0$ and $f_n(1)=1$, $f_n'(x)=\frac{k_1...k_n}{m_1...m_n}$\footnote{On each interval in $A_n$, hence of length $\frac{1}{k_1...k_n}$, $f_n$ increases by $\frac{1}{m_1...m_n}$, so indeed $f_n(1)=1$.} (we can also choose $f_n(0)=1$ and $f_n(1)=0$,$f_n'(x)=-\frac{k_1...k_n}{m_1...m_n}$). Then the sequence  $(f_n)_{n=1}^{+\infty}$ converges uniformly to some continuous, singular non-decreasing function $f_s:[0,1]\to \R$. Of course, $f_s'=0$ everywhere except on $\mathcal{C}$ (similarly, by considering $f_n(0)=0$, $f_n(1)=1$ and $f_n'(x)=-\frac{k_1...k_n}{m_1...m_n}$, we get a non-increasing $f_s$).

\begin{remark}\label{cant} {\em Let us now fix some $i_0\in\N$. For each $n>i_0$ we will define an equivalence relation $\backsim$ on the boundary $\partial A_n$ of $A_n$. Namely for $x,y\in \partial A_n$ we set
$$x\backsim y \Leftrightarrow |x-y|=\frac{r_{xy}}{k_1...k_{i_0}}$$ 
for some (unique) $r_{xy}\in \Z$, $0\leq r_{xy}\leq k_1...k_{i_0}$.\\
Notice that if $x\backsim y$ and $y\backsim z$ then since $|x-z|\leq 1$, $x\backsim z$. Hence $\backsim$ is an equivalence relation. The set $A_{i_0}$ consists of $m_1...m_{i_0}$ closed intervals , say $I_1,...I_{m_1....m_{i_0}}$ each of length $\frac{1}{k_1...k_{i_0}}$ separated by the gaps each of which has length which is a multiple of $\frac{1}{k_1...k_{i_0}}$. Clearly $\partial A_n\subset A_{i_0}$. It follows that for each $i=1,...,k_1...k_{i_0}$, each coset $E$ of $\backsim$, we have $1\leq |E\cap I_i|\leq 2$ and $|E\cap I_i|=2$ iff $0\in E$. Hence, the number of cosets is equal to $|\partial A_n\cap I_1|-1$. The latter number equals $2m_{i_0+1}...m_n-1$. Moreover, the coset $E$ containing $0$ has $2m_1...m_{i_0}$ elements, all remaining cosets having $m_1...m_{i_0}$ elements.\\
\indent Assume now that $x_0,y_0$ are in the same coset of $\backsim$ and let $I\subset \T$ be an interval. Given $t_1,t_2\in\Z$ we then have 
\begin{equation}\label{gdg}  \text{if}\;\;\;\;\; x_0+t_1\alpha,\; y_0+t_2\alpha\in I\;\;\; \text{then}\;\;\;\;\;\;\;\;\;\; |t_1-t_2|\geq \frac{C_2}{(k_1...k_{i_0})^2|I|}.
 \end{equation}
Indeed, (\ref{gdg}) follows directly from Remark \ref{dst} (with $z_1=x_0$, $z_2=y_0$) because $|x_0-y_0|=\frac{r}{k_1...k_{i_0}}$ for some $r\leq k_1...k_{i_0}$.

}

\end{remark}

\section{Weak mixing}
 We will state a criterion which implies weak mixing of such flows. First we recall a lemma.
\begin{lemma}\label{uni} (\cite{Fra-Lem}, Lemma 5.2.)  Let $T:(X,\mathscr{B},\mu)\to (X,\mathscr{B},\mu)$ be an ergodic automorphism and let $A\in \mathscr{B}$. For every $\theta, \gamma, \tau>0$ there exist $N=N(\theta,\gamma,\tau)\in \N$ and $Z=Z(\theta, \gamma, \tau)\in \mathscr{B}$ with $\mu(Z)>1-\gamma$ such that for every $M,L\in\N$ with $M,L>N$ and $L/M>\tau$ we have 
$$\left|\frac{1}{L}\sum_{j=M}^{M+L}\chi_A(T^jx)-\mu(A)\right|<\theta$$
for all $x\in Z$. \hfill $\square$
\end{lemma}
We will now prove a criterion for a special flow over an ergodic rotation to be weakly mixing.
\begin{proposition}\label{mix} Assume that $T:(X,\mathscr{B},\mu)\to(X,\mathscr{B},\mu)$ is an ergodic rotation of a compact metric Abelian group $X$ (with a translation invariant metric $d$) and $f\in L^{\infty}(X,\mathscr{B},\mu)$ positive function which is bounded away from zero. Let $P'\subset\R$ be a compact set. Assume that there exists $\eta_0>0$ such that for every $0<\eta<\eta_0$, every $\epsilon>0$ and $N\in\N$ there exist $\kappa=\kappa(\epsilon,\eta)>0$, $\delta=\delta(\epsilon,N)>0$ and a subset $Z'=Z'(\epsilon,N)\in\mathscr{B}$ with $\mu(Z')>1-\epsilon$ such that if $x,x'\in Z'$ and $0<d(x,x')<\delta$, then there are $M_1=M_1(x,x'),M_2=M_2(x,x')\geq N$, $L_1=L(x,x'), L_2=L_2(x,x')\geq N$ such that $\frac{L_1}{M_1}, \frac{L_2}{M_2}\geq \kappa$ and there exists $p=p(x,x')\in P'$ such that 
$$\frac{1}{L_1}\left|\{n\in\Z\cap[M_1,M_1+L_1] :\; |f^{(n)}(x)-f^{(n)}(x')-p|<\epsilon\}\right|>1-\epsilon,$$
and 
$$\frac{1}{L_2}\left|\{n\in\Z\cap[M_2,M_2+L_2] :\; |f^{(n)}(x)-f^{(n)}(x')-p-\eta|<\epsilon\}\right|>1-\epsilon.$$
Then the special flow $\mathcal{T}^f=(T_t^f)_{t\in \R}$ is weakly mixing.\footnote{Notice that here we must assume that $T$ is aperiodic since the assumptions of Proposition \ref{mix} are satisfied for $X$ being a finite set, while $\mathcal{T}^f$ in this case is not weakly mixing (it has discrete spectrum).} Moreover, it has cocycle Ratner's property with the set $P:=(P'+[-\eta_0,\eta_0])\setminus(-\frac{\eta_0}{4},\frac{\eta_0}{4})$.
\end{proposition}
Proof: In view of \cite{Neu},  to prove weak mixing we need to show that the equation 
\begin{equation}\label{neuma}
e^{2\pi isf(x)}=\frac{\psi(Tx)}{\psi(x)}
\end{equation}
has no measurable solution $\psi:X\to \mathbb{S}^1$ for any $s\in\R\setminus \{0\}$. We proceed by contradiction assuming that such $s,\psi$ exist. Without loss of generality, we can assume that
\begin{equation}\label{bigs} s>\frac{1}{2\eta_0}.
\end{equation}
Fix $\frac{1}{3}>\zeta>0$ (one more restriction on $\zeta$ will appear later).

 By Egorov's theorem, there exists a set $A_{\zeta}\subset X$, $\mu(A_\zeta)>1-\zeta$  such that $\psi|_{A_\zeta}$ is uniformly continuous. Therefore there exists $\delta_0>0$ such that for each $x,x',y,y'\in A_\zeta$ if $d(x,y),d(x',y')<\delta_0$ then
\begin{equation}\label{psix}\left|\frac{\psi(x')}{\psi(y')}\frac{\psi(y)}{\psi(x)}-1\right|<\frac{1}{3}\end{equation}
In view of (\ref{neuma}), for any $n\in\N$ and $x,y\in A_\zeta$ we get $e^{2\pi isf^{(n)}x}=\frac{\psi(T^nx)}{\psi(x)}$ and $e^{2\pi isf^{(n)}y}=\frac{\psi(T^ny)}{\psi(y)}$, whence
\begin{equation}\label{egu} e^{2\pi is(f^{(n)}x-f^{(n)}y)}=\frac{\psi(T^nx)}{\psi(x)}\frac{\psi(y)}{\psi(T^ny)} 
\end{equation}
Set $\eta=\frac{1}{2s}$ and fix $\epsilon>0$. Then $\kappa=\kappa(\epsilon)$ is determined. We will now use Lemma \ref{uni} with the following parameters: $A=A_\zeta$, $0<\theta<\frac{1}{3}-\zeta$, $\gamma=\zeta$ and $\tau=\kappa$. It follows that there exists $N_0$ and $Z\in\mathscr{B}$ with $\mu(Z)>1-\zeta$ such that for every $M,L>N_0$, $\frac{L}{M}\geq\kappa$ we have
\begin{equation}\label{frm}
\left|\frac{1}{L}\sum_{j=M}^{M+L}\chi_{A_{\zeta}}(T^jx)-\mu(A_{\zeta})\right|<\theta
\end{equation}
 for all $x\in Z$. By (\ref{frm}) , since $\zeta<\frac{1}{3}$, for $x\in Z$ the number of $j\in[M,M+L]$ for which $T^jx\in A_\zeta$ is at least $\frac{2}{3}L$. Hence, for any $x,y\in Z$ the number of $j\in[M,M+L]$ for which simultaneously $T^jx,T^jy\in A_\zeta$ is at least $\frac{1}{6}L$. By assumption, there exist $\delta=\delta(\epsilon,N)$ and $Z'=Z'(\epsilon, N)$, $\mu(Z')>1-\epsilon$ such that for each $x,y\in Z'$, $0<d(x,y)<\delta$ there are $M_1=M_1(x,x'),M_2=M_2(x,x')\geq N_0$, $L_1=L(x,y), L_2=L_2(x,y)\geq N_0$ such that $\frac{L_1}{M_1}, \frac{L_2}{M_2}\geq \kappa$ and there exists $p=p(x,y)\in P$ such that 
$$\frac{1}{L_1}\left|\{n\in\Z\cap[M_1,M_1+L_1] :\; |f^{(n)}(x)-f^{(n)}(y)-p|<\epsilon\}\right|>1-\epsilon,$$
and 
$$\frac{1}{L_2}\left|\{n\in\Z\cap[M_2,M_2+L_2] :\; |f^{(n)}(x)-f^{(n)}(y)-p-\frac{1}{2s}|<\epsilon\}\right|>1-\epsilon.$$ 
Notice that the set $A_\zeta\cap Z\cap Z'$ has positive measure (if $\zeta$ is small enough) and therefore, since $\mu$ is continuous, there are different points in it arbitrarily close. So fix $x,y\in  A_\zeta\cap Z\cap Z'$ such that $0<d(x,y)<\min(\delta,\delta_0)$. It follows that the number of $j\in [M_1,M_1+L_1]$ such that $T^jx,T^jy\in A_\zeta$ is at least $\frac{1}{6}L_1$. Hence there exists $j_1\in [M_1,M_1+L_1]$ such that $T^{j_1}x,T^{j_1}y\in A_\zeta$ and $|f^{(j_1)}(x)-f^{(j_1)}(y)-p|<\epsilon$. By (\ref{egu})
\begin{equation}\label{equ2}
e^{2\pi is(f^{(j_1)}x-f^{(j_1)}y)}=\frac{\psi(T^{j_1}x)}{\psi(x)}\frac{\psi(y)}{\psi(T^{j_1}y)} 
\end{equation} 
In view of (\ref{psix}) and (\ref{equ2})
$$|e^{2\pi is(f^{(j_1)}x-f^{(j_1)}y)}-1|<\frac{1}{3}.$$
On the other hand $p-\epsilon<|f^{(j_1)}x-f^{(j_1)}y|<p+\epsilon$. Hence, by letting $0<\epsilon\to 0$, we obtain $|e^{2\pi isp}-1|<\frac{1}{3}$. Doing the same for $M_2,L_2$, we get $|e^{2\pi is(p+\frac{1}{2s})}-1|<\frac{1}{3}$ or $|e^{2\pi isp}(-1)-1|<\frac{1}{3}$ which is an obvious contradiction.\\
\indent Notice that for every $p\in \R$,  $\max (|p|, |p+\frac{\eta_0}{2}|)\geq \frac{\eta_0}{4}$. So the cocycle Ratner's property is satisfied with $P$ defined in the statement of the proposition (for every $x,x'\in Z$ we choose either $p$ or $p+\frac{\eta_0}{2}$ both number being in $P'+[-\eta,\eta]$).  
\hfill $\square$

\section{Stability of Ratner's property}
In this section we show that the cocycle Ratner's property in the class of special flows over an irrational rotation by $\alpha$ having bounded partial quotients and the roof function in ${\rm BV}(\T)$ is stable under small perturbations in the variation norm. We constantly assume that $\alpha$ has bounded partial qoutients and $C:=\sup_{n\geq 1} a_n+1$.
\begin{lemma}\label{fxx} Let $f\in {\rm BV}(\T)$ and $D\in \N$. Then for every $s\in \N$
$$\max_{k<Dq_s}\sup_{\|x-y\|<1/q_s}|f^{(k)}(x)-f^{(k)}(y)|\leq 2D {\rm Var}f.$$
\end{lemma}
Proof: Let us fix $s\in \N$ and $k=cq_s+d$, with $d<q_s$, $c\leq D-1$. We will prove that for any $x,y\in\T$
\begin{equation}\label{bv} \max_{j\leq q_s}\sup_{\|x-y\|<\frac{1}{q_s}}|f^{(j)}(x)-f^{(j)}(y)|<2{\rm Var}f.
\end{equation}
 Then (\ref{bv}) and Proposition \ref{koks} will give us the required result because \begin{multline}|f^{(k)}(x)-f^{(k)}(y)|\leq \sum_{i=1}^c|f^{(q_s)}(T^{(i-1)q_s}x)-f^{(q_s)}(T^{(i-1)q_s}y)|+\\
|f^{(d)}(T^{cq_s}x)-f^{(d)}(T^{cq_s}y)|\leq
 2c{\rm Var}f+ 2{\rm Var}f\leq 2D{\rm Var}f.\end{multline}
So let us prove (\ref{bv}). Consider the points $x+i\alpha$, $i=0,...,j$. Reorder them to get $0\leq x+i_1\alpha<x+i_2\alpha<...<x+i_j\alpha$. It follows that the distance between any two such points is at least $\frac{1}{2q_s}$. Indeed, for every $u,v\in \{0,...,j\}$ we have $\|(x+u\alpha)-(x+v\alpha)\|\geq \|(u-v)\alpha\|\geq \frac{1}{2q_s}$. Set $A_1:=\{x+i_k\alpha:\; k\;\text{even}\}$ and $A_2:=\{x+i_k\alpha:\; k \;\text{odd}\}$. We get that the distance between any two points in $A_1$ and $A_2$ is at least $\frac{1}{q_s}$.  It follows that the points $x+i_k\alpha,y+i_k\alpha$ with $k$ even yield a partition of $\T$ into some intervals, and all intervals $[x+i_k\alpha,y+i_k\alpha)$ (or we have to reverse the order) are members of this partition. It follows that
$${\rm Var}f\geq \sum_{r=1}^{[j/2]}|f(x+i_{2r}\alpha)-f(y+i_{2r}\alpha)|\geq \left|\sum_{r=1}^{[j/2]}f(x+i_{2r}\alpha)-f(y+i_{2r}\alpha)\right|.$$
The same holds for the set $A_2$. Finally, we get $2{\rm Var}f\geq |f^{(j)}(x)-f^{(j)}(y)|.$
\hfill $\square$
\begin{remark}\label{rewn}\em{ Notice that by the cocycle identity it follows that for each $N\geq 1$
$$\max_{k<Dq_s}\sup_{\|x-y\|<1/q_s}|(f^{(N+k)}(x)-f^{(N+k)}(y))-(f^{(N)}(x)-f^{(N)}(y))|\leq 2D {\rm Var}f.$$
Moreover, notice that (by the proof of Lemma 5.1.) the same inequality holds if $k=Dq_s$.}
\end{remark}

\begin{theorem}\label{sta}
Let $P\subset \R\setminus\{0\}$ be a nonempty compact set. Let $\frac{1}{2}>\eta>0$ be such that $P\subset \R\setminus(\eta,\eta)$ and  $(T_t^f)_{t\in\R}$ be a special flow over an irrational rotation by $\alpha$ with the roof function $f\in L^\infty(X,\mathscr{B},\mu)$. Assume that $(T_t^f)_{t\in\R}$ satisfies the cocycle Ratner's property (see Definition \ref{krp}) and assume additionally that there exist $R_1,R_2>0$ such that $R_1\leq M(x,x')\|x-x'\|\leq R_2$. Moreover, assume that for every $n\in\Z\cap[M,M+L]$, we have $|f^{(n)}(x)-f^{(n)}(x')-p|<\epsilon$ for every $x,x'\in\T$ (this does not depend on $\epsilon, N,$ see Introduction). Then the cocycle Ratner's property is stable under bounded variation perturbations, more precisely if $g\in {\rm BV}(\T)$ with ${\rm Var}g<\frac{1}{8CR_2}\eta$ then $(T_t^{f+g})_{\in \R}$ has the cocycle Ratner's property with some compact set $P'\subset \R\setminus\{0\}$.
\end{theorem}
Proof: Fix $\epsilon>0$ and $N\in \N$. Let us define $\kappa_{f+g}=\kappa_{f+g}(\epsilon):=\frac{\epsilon}{4}\min(\kappa_f(\epsilon/2),\frac{1}{CR_2})$ and $\delta_{f+g}=\delta_{f+g}(\epsilon,N)=\delta_f(\epsilon/2,N)$ (cf.  footnote~\ref{zale}). Let us take $0<\|x-x'\|<\delta_{f+g}$. By assumption, there exist $M_f(x,x'), L_f(x,x')$ and $p_f(x,x')$ such that $M_f>N$, $1>\frac{L_f}{M_f}>\kappa_f$ and for every $k\in [M_f, M_f+L_f]$
$$|f^{(k)}(x)-f^{(k)}(x')-p_f|<\epsilon/2.$$
Let $s\in \N$ be a unique natural number such that $\frac{1}{q_{s+1}}\leq\|x-x'\|<\min(\frac{1}{q_s}\delta_{f+g})$. Let us now denote $L'_f=\min(q_s, L_f)$  consider the interval $[M_f,M_f+L'_f]$. It follows by the proof of (\ref{bv}) in Lemma \ref{fxx} ($L'_f\leq q_s$, $\|T^{M_f}x-T^{M_f}x'\|=\|x-x'\|<\frac{1}{q_s}$) that
\begin{multline}\label{prze} \frac{1}{4}>2 {\rm Var}g>|g(T^{M_f}x)-g(T^{M_f}x')|+|g(T^{M_f+1}x)-g(T^{M_f+1}x')|+\\
...+|g(T^{M_f+L'_f}x)-g(T^{M_f+L'_f}x')|
\end{multline}
Therefore, there exist an interval $[M_f,M_f+L'_f]\supset I_g=[a_g,b_g]$ of length $[\epsilon L'_f]$ such that 
$$\frac{\epsilon}{4}>|g(T^{a_g}x)-g(T^{a_g}x')|+|g(T^{a_g+1}x)-g(T^{a_g+1}x')|+...+|g(T^{b_g}x)-g(T^{b_g}x')|;$$
in particular, we get that for every $k\in I_g$
\begin{equation}\label{smalg} |g^{(k)}(x)-g^{(k)}(x')-g^{(a_g)}(x)-g^{(a_g)}(x')|<\frac{\epsilon}{2}.
\end{equation}

Let $M_{f+g}=a_g, L_{f+g}=b_g-a_g\geq[\epsilon L'_f]\geq \frac{\epsilon}{2}L'_f=\frac{\epsilon}{2}\min(L_f,q_s)$. Then, by definition, $N\leq M_f\leq M_{f+g}\leq M_f+L_f<2M_f$, $\frac{L_{f+g}}{M_{f+g}}\geq \frac{\epsilon}{2} \frac{\min (L_f,q_s)}{2M_f}\geq \frac{\epsilon}{4}\min(\frac{L_f}{M_f},\frac{q_s}{M_f})=\kappa_{f+g}$. Moreover, let $P':=P+[-\eta/2,\eta/2]\subset \R\setminus (-\eta/2,\eta/2)$ and $p_0(x,x'):=g^{(M_{f+g})}(x)-g^{(M_{f+g})}(x')$. It follows by Lemma \ref{fxx} and the fact that ${\rm Var}g<\frac{\eta}{8CR_2}$ ($M_{f+g}\leq 2M_f\leq 2CR_2q_s$) that $|p_0|\leq 4CR_2 \frac{\eta}{8CR_2}=\eta/2$. Let $p'(x,x'):=p_f(x,x')+p_0(x,x')\in P'$. It view of (\ref{smalg}), it follows that for every $k\in [M_{f+g}, M_{f+g}+L_{f+g}]\subset[M_f,M_f+L_f]$
\begin{multline}
|(f+g)^{(k)}(x)-(f+g)^{(k)}(x')-p'|\leq|f^{(k)}(x)-f^{(k)}(x')-p_f|+|g^{(k)}(x)-g^{(k)}(x')-p_0| <\\
\epsilon/2+\epsilon/2=\epsilon.
\end{multline}
This completes the proof. \hfill $\square$

\section{Weak mixing and Ratner's property for some roof functions in ${\rm BV}(\T)$}

In this section we prove weak mixing and Ratner's property of a special flow $\mathcal{T}^f$, where $f\in {\rm BV}(\T)$ is of the form $f=f_a+f_j+f_s+S'\{\cdot\}$, with $f_a$ is absolutely continuous on $\T$, $f_j$ is piecewise constant (finitely many jumps) and $f_s$ is a singular, continuous on $[0,1]$, function which comes from a quasi-similar Cantor set ($f_s=f(\mathcal{C})$). Consider $\mathcal{C}=\cap_{i=1}^{+\infty}A_i$ (with bounded sequences $(m_i)_{i\geq 1}$, $(k_i)_{i\geq 1}$, see Introduction). Take any two points $x<y\in \T$ and let $n=n(x,y)\in\N$ be a unique natural number such that 
\begin{equation}\label{star}\frac{1}{k_1...k_{n+1}}<\|x-y\|\leq \frac{1}{k_1...k_n},\end{equation} then $f_s(x)\neq f_s(y)$ iff $[x,y]\cap \partial A_{n+1}\neq \emptyset$. Indeed, at step $n\geq 1$ of the construction of $\mathcal{C}$, $[0,1]$ is divided into $k_1...k_{n+1}$ intervals of equal length $\frac{1}{k_1...k_{n+1}}$, $m_1...m_{n+1}$ of which, say $I_1,..., I_{m_1...m_{n+1}}$, yield $A_{n+1}$. Now, $x$ and $y$ cannot belong to the same $I_j$, hence either
\begin{itemize}
    \item $[x,y]\cap\partial A_{n+1}=\emptyset$, i.e.  $[x,y]\cap I_j=\emptyset$ for each $j$ and then $f_s(x)=f_s(y)$
\end{itemize}
or
\begin{itemize}
	 \item  $[x,y]\cap\partial A_{n+1}\neq\emptyset$ i.e. there is some interval $I_{i_0}$ such that $[x,y]\cap I_{i_0}\neq \emptyset$. In this case $f_s(x)\neq f_s(y)$.
\end{itemize}
Moreover, it follows from (\ref{star}) that $|f_s(y)-f_s(x)|\leq \frac{1}{m_1m_2...m_n}$ (unless $0\in [x,y]$). Indeed, by definition, for every $n\in \N$, $f_s$ is constant on each connected component of $\mathcal{C}\setminus A_n$, and for every chosen, at step $n$, interval $I'=[a',b')$, $f_s(b')-f_s(a')=\frac{1}{m_1...m_n}$.  

 To prove both weak mixing and WR-property for $(T_t^f)$ we need a lemma.
\begin{lemma}\label{abs} {\em(\cite{Fr-Lem}, Lemma 6.1.)}
 Let $T:\mathbb{T}\to \mathbb{T}$ be the rotation by an irrational $\alpha$ with bounded partial quotients and let $f:\mathbb{T}\to \R$ be absolutely continuous with zero mean. Then 
$$\sup_{0\leq n<q_{s+1}}\sup_{\|y-x\|<\frac{1}{q_s}}|f^{(n)}(y)-f^{(n)}(x)|\to 0,\;\; \mbox{as}\;s\to +\infty.$$
\end{lemma}
\begin{remark}\label{abs2} {\em Following (step by step) the proof of Lemma 6.1. in \cite{Fr-Lem} one can prove that for every $b\in \N$ we have
$$\sup_{0\leq n<q_{s+b}}\sup_{\|y-x\|<\frac{1}{q_s}}|f^{(n)}(y)-f^{(n)}(x)|\to 0,\;\; \mbox{as}\;s\to +\infty.$$ }
\end{remark}
Let us denote the (finite) set of discontinuities of $f$ by $\{\beta_i\}_{i=0}^k$ with the corresponding set of jumps $\{d_i\}_{i=0}^k$ (we recall that $d_0$ may be equal to $0$).\\
\indent\textbf{Proof of Theorem \ref{wmrp}}:

We will use Proposition \ref{mix} with $P:=[-2C{\rm Var} f, 2C{\rm Var}f]$. It follows that $f=f_{ac}+f_j+f_s +\tilde{S}\{\cdot\}+c$, where $f_{ac}$ is absolutely continuous with zero mean $(c:=\int_\T f_a(x)d\mu)$. Let $f_{pl}(x):=f_j(x)+\tilde{S}\{x\}+c$, so 
\begin{equation}\label{f-} f=f_{ac}+f_{s}+f_{pl}.
\end{equation} 
 Assume that $\tilde{S}>0$. Fix $0<\eta\leq \min(\frac{\tilde{S}}{5C(2C+1)(k+1)},\frac{1}{8}), \epsilon>0$ and $N\geq 1$. By Lemma \ref{abs}, there exists $s_0$ such that for every $s\geq s_0$, we have  
\begin{equation}\label{f-ac}\sup_{0\leq n<q_{s+1}}\sup_{\|y-x\|<\frac{1}{q_s}}|f_{ac}^{(n)}(y)-f_{ac}^{(n)}(x)|<\frac{\epsilon}{8}.\end{equation}
 Let $n_0=n_0(\epsilon)$ be a unique natural number such that $\frac{2K}{m_1...m_{n_0}}\leq\epsilon/4<\frac{2K}{m_1...m_{n_0+1}}$\footnote{$K$ is a constant such that $m_i\leq K$ for every $i\in \N$.}.
ad $1.$ First we have to define $\kappa=\kappa(\epsilon)$\footnote{Note that here $\kappa$ does not depend on $\eta>0$.}. Set \\
$\kappa(\epsilon):=\min(\frac{C_2}{C(k_1...k_{n_0})^2}, \frac{\epsilon}{4\tilde{S}C},\frac{1}{2(2C+1)(k+1)}$) (see Lemma \ref{dist}). Let $s_1\in \N$ be the smallest natural number with $q_{s_1}\geq \max(\tilde{S}q_{s_0},\frac{16\tilde{S}}{\epsilon}, \frac{N}{\kappa})$ and
let $\delta(\epsilon,N):=\min(\frac{1}{q_{s_1}},\frac{1}{k_1...k_{n_0}})$. Take any $x,y\in \T$, $\|x-y\|<\delta$. Let $s\in\N$ be a unique natural number such that $\frac{1}{q_{s+1}}\leq\|x-y\|<\frac{1}{q_s}$. Consider now the points $\T\ni \beta_i-j\alpha$, $i=0,...,k$, $j=0,...,q_{s+1}-1$. It follows that for every $i$, the number of $t=0,...,q_{s+1}-1$ such that $\beta_{i}-t\alpha\in[x,y)$ is at most $2C+1$. So we can divide the time interval $[q_s,q_{s+1}]$ into $(2C+1)(k+1)$ clopen intervals $I_1,...,I_{(2C+1)(k+1)}$ such that for every $u=1,...,(2C+1)(k+1)$, $v\in {\rm int} \;I_u$ and $i=0...,k$, we have $\beta_i-v\alpha \notin [x,y)$. One of them, say $I=[M_0, K_0]$, has length at least $\frac{1}{(2C+1)(k+1)}(q_{s+1}-q_s)$. It follows that $M_0<q_{s+1}$, let 
\begin{equation}\label{kap} L_0:=\kappa M_0<\kappa q_{s+1}\leq \kappa Cq_s\leq \frac{C_2q_s}{(k_1...k_{n_0})^2}.
\end{equation}
Moreover, $M_0>q_s>N$, $L_0\geq \kappa M_0>N$. We define $p(x,y):=f^{(M_0)}(y)- f^{(M_0)}(x)$, then $|p|<2C\!{\rm Var}f$ (by Lemma \ref{fxx}). Denote $n:=n_{xy}$ to be a unique natural number such that $\frac{1}{k_1...k_{n}}\leq\|x-y\|<\frac{1}{k_1...k_{n-1}}$. If  $f_s(x)\neq f_s(y)$ then $ \partial A_n\cap [x,y]\neq \emptyset$. Recall that $|\partial A_n|=2m_1...m_n$. Take $z\in \partial A_n$. It follows by Remark \ref{dst}, with $z_1=z_2, r=0, q=1$ and $I=[x,y]$ and by (\ref{kap}) that there exists at most one $[M_0,M_0+L_0]\ni v_0$ such that $z+v_0\alpha\in [x,y]$. Note that for every $v\in [M_0,M_0+L_0]$ we have  $|f_s(y+v\alpha)-f_s(x+v\alpha)|<\frac{1}{m_1...m_{n-1}}$ (because $-\{i\alpha\}\notin [x,y]$ for $i\in [M_0,M_0+L_0]$). Consider now Remark \ref{cant} with $i_0=n_0$ and the corresponding equivalence relation $\backsim$ for $n$ (note that $n\geq n_0$). It follows by (\ref{kap}) and (\ref{gdg}) that there is at most one point $z_i\in \partial A_n$ in each coset for which there exists $w_i\in [M_0,M_0+L_0]$ such that $z_i+w_i\alpha\in [x,y]$. So there are at most $2m_{n_0+1}...m_n-1$ points $z\in\partial A_n$ for which there exists a (unique) $w_z\in [M_0,M_0+L_0]$ such that $z+w_z\alpha \in [x,y]$. For each such $w_z\in [M_0,M_0+L_0]$, we have $0 \neq f_s(y+w_z\alpha)-f_s(x+w_z\alpha)<\frac{1}{m_1...m_{n-1}}$.
There are at most $2m_{n_0+1}...m_n-1<2m_{n_0+1}...m_n$ such $w=w_z\in [M_0,M_0+L_0]$. Hence, for every $r\in[M_0,M_0+L_0]$ (remembering that $f_s(T^{M_0+j}y)=f_s(T^{M_0+j}x)$ for $M_0+j\neq w_z$ for each $z$ above)
\begin{multline}\label{f-s}|f_s^{(r)}(y)- f_s^{(r)}(x)-(f_s^{(M_0)}(y)- f_s^{(M_0)}(x))|=|f_s^{(r-M_0)}(T^{M_0}y)-f_s^{(r-M_0)}(T^{M_0}x)|\leq\\
\sum_{j=0}^{r-M_0-1}|f_s(T^{M_0+j}y)-f_s(T^{M_0+j}x)|\leq \sum_{j=0}^{L_0-1}|f_s(T^{M_0+j}y)-f_s(T^{M_0+j}x)|\leq\\
2m_{n_0+1}...m_n\frac{1}{m_1...m_{n-1}}=
 \frac{2m_n}{m_1...m_{n_0}}\leq\frac{\epsilon}{4}.
\end{multline}
\begin{remark}\label{potr} { \em Note that if $I=[A,B]$ is an interval such that $|I|\leq \frac{C_2q_s}{(k_1...k_{n_0})}$ and for any $v\in [A,B]$ we have $v\alpha \notin [x,y)$ then for every $r\in [A,B]$ 
$$|f_s^{(r)}(y)- f_s^{(r)}(x)-(f_s^{(A)}(y)- f_s^{(A)}(x))|\leq\frac{\epsilon}{4}.$$
 Indeed, these are the only assumptions to prove (\ref{f-s}).}
\end{remark}

Moreover, it is easy to see (because $f_{pl}$ is piecewise linear and there is no discontinuity point of $f$ in $[T^rx,T^ry]$ for $r\in [M_0,M_0+L_0]$) that for $r\in [M_0,M_0+L_0]$
\begin{multline}\label{f-pl}|(f_{pl}^{(r)}(y)- f_{pl}^{(r)}(x))-(f_{pl}^{(M_0)}(y)- f_{pl}^{(M_0)}(x))| \leq \sum_{j=0}^{r-M_0-1}(f_{pl}(T^{M_0+j}(y))-f_{pl}(T^{M_0+j}(x)))\leq\\
 (r-M_0) \tilde{S}\|y-x\|<L_0\tilde{S}\|y-x\|\leq\tilde{S} \kappa M_0\frac{1}{q_s}\leq \kappa\tilde{S} C<\frac{\epsilon}{4}.\end{multline}

Finally, by (\ref{f-}), (\ref{f-ac}), (\ref{f-s}), (\ref{f-pl}), we get for $r\in [M_0,M_0+L_0]$

\begin{multline}\label{sds}|f^{(r)}(y)- f^{(r)}(x)-p|=|(f^{(r)}(y)- f^{(r)}(x))-(f^{(M_0)}(y)- f^{(M_0)}(x))|=\\
|(f_{ac}^{(r)}(y)- f_{ac}^{(r)}(x))-(f_{ac}^{(M_0)}(y)- f_{ac}^{(M_0)}(x))+(f_s^{(r)}(y)- f_s^{(r)}(x))-(f_s^{(M_0)}(y)- f_s^{(M_0)}(x))+\\
(f_{pl}^{(r)}(y)- f_{pl}^{(r)}(x))-(f_{pl}^{(M_0)}(y)- f_{pl}^{(M_0)}(x))|\leq\\
\frac{\epsilon}{8}+\frac{\epsilon}{4}+\frac{\epsilon}{4}<\frac{5}{8}\epsilon<\epsilon.\end{multline}

\begin{remark}\label{emm} {\em Notice that if $M'<q_{s+1}$ is such that for $v'\in[M',M'+\kappa M']$ we have $\beta_i+v'\alpha \notin [x,y)$, $i=0,...,k$ and $p':=f^{(M')}(y)- f^{(M')}(x)$, then for any $r'\in [M',M'+\kappa M']$ $|f^{(r')}(y)- f^{(r')}(x)-p'|\leq\frac{5}{8}\epsilon$. Indeed, these are the only assumptions on $M_0$ to prove (\ref{sds}).
}
\end{remark}
Consider again the  interval $I=[M_0,K_0]$. We recall that none of the points of the form $\beta_i-v\alpha$, $i=0,...,k$ and $v\in [M_0,K_0]$ belongs to $[x,y]$. Therefore, using the fact that $\tilde{S}>0$, $f_s$ non-decreasing, for every $r\in [M_0, K_0]$, we have
\begin{equation}\label{cmp}\tilde{S}\|y-x\|\leq (f_{pl}+f_s)(T^ry)-(f_{pl}+f_s)(T^rx)\leq \tilde{S}\|y-x\|+ \frac{1}{m_1...m_n}\leq \frac{\epsilon}{8}.\end{equation}
The length of $I$ is at least $\frac{1}{(2C+1)(k+1)}(q_{s+1}-q_s)$, so by (\ref{f-ac}) and (\ref{cmp}) for every $\ell> M_0$, we have 
   \begin{multline}\label{elll}|(f^{(\ell)}(y)- f^{(\ell)}(x))-(f^{(M_0)}(y)- f^{(M_0)}(x))|=\\
\left|(f^{(\ell-M_0)}_{ac}(T^{M_0}y)-f^{(\ell-M_0)}_{ac}(T^{M_0}x))+((f_{pl}+f_s)^{(\ell-M_0)}(T^{M_0}y)- (f_{pl}+f_s)^{(\ell-M_0)}(T^{M_0}x))\right|\geq\\
\left|-\frac{\epsilon}{4} +\sum_{i=M_0}^{\ell}(f_{pl}+f_s)(T^{M_0+i}y)-(f_{pl}+f_s)(T^{M_0+i}x)\right|\geq\\
(\ell-M_0) \tilde{S}\|y-x\|-\frac{\epsilon}{4}.
\end{multline}
 In particular, for $\ell_0=M_0+ [\frac{K_0}{2}]$, $|(f^{(\ell_0)}(y)- f^{(\ell_0)}(x))-(f^{(M_0)}(y)- f^{(M_0)}(x))|>\frac{\tilde{S}}{4C(2C+1)(k+1)}\geq\eta$.\\
\indent Notice that if we set $H(r):=(f^{(r)}(y)- f^{(r)}(x))-(f^{(M_0)}(y)- f^{(M_0)}(x))$ then for every $r\in[ M_0,M_0+\ell_0]$
\begin{equation}\label{ha} |H(r+1)-H(r)|<\frac{\epsilon}{4}.
\end{equation}
Indeed, $|H(r+1)-H(r)|=|f(T^ry)-f(T^rx)|$ and (\ref{ha}) follows by (\ref{cmp}) and (\ref{f-ac}). It follows by (\ref{elll}) and (\ref{ha}) (by considering $r=M_0,...,M_0+\ell_0$) that there exists $R_0\in [M_0,\ell_0]$ such that  
$$\eta+\frac{\epsilon}{4} \geq |(f^{(R_0)}(y)- f^{(R_0)}(x))-(f^{(M_0)}(y)- f^{(M_0)}(x))|\geq \eta.$$
Define $M_1:=R_0$ and $L_1:=\kappa M_1$. Then $M_1>q_s>N$, $L_1\geq \kappa M_1>N$.
Hence 
\begin{equation}\label{mimi}|f^{(M_1)}(x)- f^{(M_1)}(y)-p-\eta|<\frac{\epsilon}{4}.\end{equation}

It follows by Remark \ref{emm} and (\ref{mimi}) that for any $r\in[M_1,M_1+L_1]$
\begin{multline}\label{tre}|f^{(r)}(x)- f^{(r)}(y)-p-\eta|\leq |(f^{(r)}(x)- f^{(r)}(y))-(f^{(M_1)}(x)- f^{(M_1)}(y))|+\\
|f^{(M_1)}(x)-f^{(M_1)}(y)-p-\eta|<\frac{5}{8}\epsilon+\frac{1}{4}\epsilon<\epsilon.\end{multline}

ad 2. First note that there exist $b\in \N$, b<C, such that for every $s\geq 1$, $\left[\frac{q_{s+b}}{q_{s+1}}\right]\geq 6C+4$. By Remark \ref{abs2}, there exists $s_0$ such tht for every $s\geq s_0$ we have 
\begin{equation}\label{f-ac2}\sup_{0\leq n<q_{s+b}}\sup_{\|y-x\|<\frac{1}{q_s}}|f_{ac}^{(n)}(y)-f_{ac}^{(n)}(x)|<\frac{\epsilon}{8}.\end{equation}
 
Define
$\kappa=\kappa(\epsilon):=\min(\frac{C_2}{C^b(k_1...k_{n_0})^2},\frac{1}{(2C+1)^{b}})$.
Let $s_1\in \N$ be the smallest natural number with $q_{s_1}\geq \max(q_{s_0},N)$ and
let $\delta(\epsilon,N):=\frac{1}{q_{s_1}}$. Take any $x,y\in \T$, $\|x-y\|<\delta$. Let $s\in\N$ be a unique natural number such that $\frac{1}{q_{s+1}}\leq\|x-y\|<\frac{1}{q_s}$. Consider the time interval $[q_s,q_{s+b}]$. Denote by $(q_s\leq) R_0<...<R_t(:=q_{s+b}-j_t)$ all natural numbers in $[q_{s},q_{s+b}]$ for which $-R_i\alpha\in[x,y)$. These numbers divide $[q_s,q_s+b]$ into some clopen subintervals $I_0,...,I_t$, $C q_s\geq q_{s+1}\geq |I_i|\geq C_2 q_s$ (by Remark \ref{dst}); moreover ($6C+4\leq [\frac{q_{s+b}}{q_{s+1}}]\leq t\leq (2C+1)^b$). In particular, $j_t\leq q_{s+1}$, $R_0\leq 2q_{s+1}$. It follows that $|(f_s^{(R_i+1)}(x)- f_s^{(R_i+1)}(y))-(f_s^{(R_i)}(x)- f_s^{(R_i)}(y))|>\frac{7}{8},$
 for $i=0,...,t$ (because $-R_i\alpha\in [x,y]$). Moreover, there exists an $i\in\{0,...,t-1\}$ such that $(f_s^{(R_{i+1})}(x)- f_s^{(R_{i+1})}(y))-(f_s^{(R_i)}(x)- f_s^{(R_i)}(y))\geq -\frac{3}{4}$ (then  $(f_s^{(R_{i+1})}(x)- f_s^{(R_{i+1})}(y))-(f_s^{(R_i+1)}(x)- f_s^{(R_i+1)}(y))\geq \frac{7}{8}-\frac{3}{4}= \frac{1}{8}$ ). Indeed, by Proposition \ref{koks} and Lemma \ref{fxx} we have 
\begin{multline} 2C+2=(2C+2){\rm Var}f_s\geq  |f_s^{(q_{s+b})}(x)- f_s^{(q_{s+b})}(y)|+|f_s^{(-j_t)}(T^{q_{s+b}}x)- f_s^{(-j_t)}(T^{q_{s+b}}y)\footnote{$f_s^{(-j_t)}(T^{q_{s+b}}x)=-f_s^{(j_t)}(T^{q_{s+b}-j_t}x)$}|\geq\\ |f_s^{(q_{s+b}-j_t)}(x)- f^{(q_{s+b}-j_t)}(y)|=|f^{(R_{t})}(x)- f^{(R_{t})}(y)|=\\
\big|\sum_{w=0}^{t-1}(f_s^{(R_{w+1})}(x)- f_s^{(R_{w+1})}(y))-(f_s^{(R_w)}(x)- f_s^{(R_w)}(y))+\\
 (f_s^{(R_0)}(x)- f_s^{(R_0)}(y))\big|> t|-\frac{3}{4}|-2C>2C+2. 
\end{multline}
Let us define $M_0:=R_i+1<q_{s+b}, L_0:=\kappa M_0<\kappa q_{s+b}\leq \frac{C_2q_s}{(k_1...k_{n_0})^2}$. Then $M_0>q_s>N$, $L_0\geq \kappa M_0>N$. It follows by Remark \ref{potr} that for every $r\in [M_0,M_0+L_0]$ 
$$|(f_s^{(r)}(x)- f_s^{(r)}(y))-(f_s^{(M_0)}(x)- f_s^{(M_0)}(y))|<\frac{\epsilon}{4}.$$
Let us define $p=p(x,y):= f^{(M_0)}(x)- f^{(M_0)}(y)$. Then for $r\in[M_0,M_0+L_0]$ we have (by (\ref{f-ac2}))
$$|(f^{(r)}(x)- f^{(r)}(y))-(f^{(M_0)}(x)-f^{(M_0)}(y))|\leq \frac{\epsilon}{2}.$$
Set $H(r):= (f^{(r)}(x)- f^{(r)}(y))-(f^{(M_0)}(x)-f^{(M_0)}(y))$ for $r\in [M_0,R_{i+1}]$ ($H(R_i+2)=0, H(R_{i+1})\geq \frac{1}{4}\geq 2\eta$). It follows by the choice of $i$, (\ref{cmp}) (with $f_{pl}=0$, $\tilde{S}=0$) that there exists $R_\eta\in [M_0,R_{i+1}]$ such that 
$|f_s^{(R_\eta)}(x)- f_s^{(R_\eta)}(y)-p-\eta|<\frac{\epsilon}{2}$. Then defining $M_1:=R_\eta$, $L_1:=\kappa M_1$ ($M_1>q_s>N$, $L_1\geq \kappa M_1>N$), we get (proceeding like in (\ref{tre}) and by (\ref{f-ac2})) that for every $r\in [M_1,M_1+L_1]$
$$|f^{(r)}(x)- f^{(r)}(y)-p-\eta|<\epsilon.$$
The proof of Theorem \ref{wmrp} is complete. \hfill $\square$

\section{Absence of partial rigidity and mild mixing.}
In this section, we will show the absence of partial rigidity of some special flows over irrational rotation by $\alpha$ having bounded partial qoutients ($\sup_{n\geq 1}a_n+1<+\infty$) and roof functions of the form $f=f_a+f_j+f_s+\tilde{S}\{\cdot\}$, where $f_j$ has finitely many jumps,  $f_s=f(\mathcal{C})$ is non-decreasing and $\tilde{S}>0$~\footnote{The case $\tilde{S}<0$ and $f_s$ decreasing goes analogously.}.
The proof of the theorem below can be obtained by a repetition of the proof of Theorem  7.1. in \cite{Fr-Lem}. 
\begin{theorem}\label{rig} Assume $T:\T\to\T$ is an ergodic rotation by $\alpha$ having bounded partial quotients. Suppose $f$ is of the above form. Then the special flow $(T_t^f)_{t\in\R}$ is not partially rigid.
\end{theorem}
Combining Theorems \ref{rig} and \ref{wmrp} with the main argument used in \cite{Fr-Lem} to prove Theorem 7.2 therein, we conclude in the following result.
\begin{corollary} 
Suppose that $T:\T\to \T$ is the rotation by an irrational number $\alpha$
with bounded partial quotients and $f:\T\to\R$ is a positive, bounded away from zero function of the form $f:=f_a+f_s+f_j+\tilde{S}\{\cdot\}$ with $\tilde{S}\neq 0$ and $f_s=f_s(\mathcal{C})$ for some quasi-similar Cantor set $\mathcal{C}$. Then $(T_t^f)$ is mildly mixing. 
\end{corollary}

Institute of Math.\\
Polish Academy of Sciences,\\
Śniadeckich 8,\\
00-950 Warszawa, Poland\\
adkanigowski@gmail.com

\end{document}